\documentclass{siamltex}

\usepackage{amsmath,amssymb}
\usepackage{graphicx}
\usepackage{hyperref}
\usepackage{booktabs}
\usepackage{subcaption}
\usepackage{tikz}
\usetikzlibrary{patterns,arrows.meta,decorations.pathreplacing}
\usepackage{algorithm}
\usepackage{algpseudocode}
\usepackage{xcolor}

\newcommand{\dd}{\mathrm{d}}
\newcommand{\RR}{\mathbb{R}}

\begin{document}

\title{Semi-Lagrangian Discontinuous Galerkin Method with Adaptive Mesh
       Refinement for the Vlasov--Poisson System in 1X+3V}

\author{Mark~F.~Adams\thanks{Computational Research Division,
  Lawrence Berkeley National Laboratory, Berkeley, CA 94720, USA
  (\texttt{mfadams@lbl.gov}).}}

\maketitle

\begin{abstract}
We extend the semi-Lagrangian discontinuous Galerkin (SLDG) method of
Einkemmer to velocity grids with adaptive mesh refinement (AMR) and to
three-dimensional velocity space.  The original SLDG formulation assumes
uniform cell widths, which permits the overlap matrices to be precomputed
once per fractional shift and reused for every cell.  On an adaptively
refined mesh, neighboring cells may differ in size, invalidating this
assumption.  We develop a hybrid sweep strategy: conforming cells in the
mesh interior use precomputed per-level overlap matrices (the fast path),
while nonconforming cells at refinement boundaries evaluate generalized
overlap integrals on the fly (the slow path).  A compressed sparse row
(CSR) pencil data structure organizes the dimensional splitting along each
velocity coordinate, with weighted accumulation for coarse cells that
appear in multiple pencils.  The method is extended from one to three
velocity dimensions using tensor-product DG elements on hexahedral cells
provided by PETSc's PetscFE class.  We verify the solver on
the standard Landau damping benchmark in 1X+3V, demonstrating correct
damping rates, exact mass conservation, and convergence behavior with
polynomial degree and AMR refinement level.
\end{abstract}


\section{Introduction}
\label{sec:intro}

Kinetic simulation of magnetized plasmas requires solving the Vlasov
equation in phase space, coupling spatial transport with velocity-space
dynamics including electromagnetic forces and Coulomb collisions.  The full
Vlasov--Maxwell--Landau system evolves the distribution function
$f_\alpha(t,\vec{x},\vec{v})$ for each species~$\alpha$ according to
\begin{equation}\label{eq:vml}
  \frac{\partial f_\alpha}{\partial t}
  + \vec{v} \cdot \nabla_x f_\alpha
  + \frac{q_\alpha}{m_\alpha}(\vec{E} + \vec{v}\times\vec{B})
    \cdot \nabla_v f_\alpha
  = \sum_\beta C[f_\alpha, f_\beta],
\end{equation}
where $C[f_\alpha,f_\beta]$ is the Landau collision operator
\cite{Landau1936,HirvijokiAdams2017}.  Operator splitting separates the
collisionless Vlasov--Maxwell transport from the velocity-space collision
step, which can be advanced independently at each spatial point.  A
structure-preserving Landau collision operator with adaptive mesh refinement
(AMR) and batched linear solvers is available in PETSc
\cite{AdamsLandau2025,AdamsHirvijoki2017,AdamsIPDPS2022}.

The collisionless transport step requires solving the advection equation in
phase space.  Semi-Lagrangian methods are widely used for this purpose
because they are unconditionally stable with respect to the CFL condition,
allowing large time steps \cite{Sonnendrucker1999,Sonnendrucker2011,ChengKnorr1976}.
Traditional semi-Lagrangian schemes based on polynomial or spline
interpolation suffer from either large communication stencils or global
solves that limit parallel scalability \cite{Einkemmer2016}.

The semi-Lagrangian discontinuous Galerkin (SLDG) method, developed by
Crouseilles, Mehrenberger, and Vecil \cite{CrouseillesMehrenbergerVecil2011}
and analyzed by Einkemmer and Ostermann \cite{EinkemmerOstermann2014},
offers an attractive alternative.  The method represents the distribution
function as a piecewise polynomial on each cell, translates it along
characteristics, and projects back onto the DG basis via overlap integrals.
The resulting scheme is conservative by construction, requires data from at
most two adjacent cells regardless of polynomial degree, and exhibits low
numerical diffusion \cite{EinkemmerOstermann2014}.
Einkemmer \cite{Einkemmer2016} demonstrated excellent parallel scalability
of this approach on distributed-memory systems for uniform Cartesian grids.

However, the original SLDG formulation assumes uniform grid spacing, which
allows the overlap matrices $A^\alpha$ and $B^\alpha$ to be precomputed
once per fractional shift~$\alpha$ and reused for all cells.  This
assumption breaks down on adaptively refined meshes, where neighboring
cells may have different sizes.  Extending the SLDG method to AMR grids is
the primary contribution of this paper.

\paragraph{Contributions.}
This paper extends the SLDG method of Einkemmer \cite{Einkemmer2016} in
two directions:
\begin{enumerate}
\item \textbf{AMR velocity space.}  We develop a hybrid SLDG sweep that
  handles adaptively refined velocity grids built on PETSc's p8est
  infrastructure \cite{BursteddeWilcoxGhattas2011}.  Conforming cells in
  the mesh interior use precomputed per-level overlap matrices (the fast
  path), while nonconforming cells at refinement boundaries use generalized
  overlap integrals computed on the fly (the slow path).  A compressed
  sparse row (CSR) pencil data structure organizes the sweep along each
  velocity dimension, with weighted accumulation for coarse cells that
  appear in multiple pencils.

\item \textbf{3V velocity space.}  We extend the 1X+1V formulation to
  1X+3V with a full three-dimensional Cartesian velocity space
  $(v_x, v_y, v_z)$, using tensor-product DG elements and dimensional
  splitting along pencils extracted from the 3D mesh topology.
\end{enumerate}
The solver is implemented in PETSc \cite{PETSc2024} and is designed to
couple with PETSc's existing Landau collision operator
\cite{AdamsLandau2025} for full Vlasov--Maxwell--Landau simulations.  We
verify the implementation on the standard Landau damping benchmark,
demonstrating correct damping rates, exact mass conservation, and the
expected convergence behavior with polynomial degree and mesh refinement.

\paragraph{Scope and roadmap.}
This paper focuses on the collisionless transport step of the
Vlasov--Maxwell--Landau system: the Vlasov--Poisson equations in 1X+3V
(one spatial dimension, three velocity dimensions).  The full system
\eqref{eq:vml} includes the Lorentz force from electromagnetic fields and
the Landau collision operator on the right-hand side.  The Strang splitting
framework (Section~\ref{sec:strang}) cleanly separates these components:
the collision operator is already available in PETSc
\cite{AdamsLandau2025}, and the electromagnetic field solver is a natural
extension of the Poisson solver presented here.  By developing and
validating the SLDG transport solver independently on the Vlasov--Poisson
benchmark, we establish the accuracy and conservation properties of the
advection scheme before coupling it with the collision and field solvers in
future work.

\paragraph{Related work.}
The semi-Lagrangian approach to the Vlasov equation was established by
Sonnendr{\"u}cker et al.\ \cite{Sonnendrucker1999}, who introduced the
splitting of the six-dimensional phase-space advection into a sequence of
lower-dimensional sweeps and demonstrated that cubic spline interpolation
along characteristics yields an efficient and accurate scheme for the
Vlasov--Poisson system.
Adaptive mesh refinement for Vlasov solvers has been explored in several
contexts.  Besse and Segr{\'e} \cite{BesseSegre2008} developed an
adaptive semi-Lagrangian method on locally refined grids using cubic
spline interpolation, demonstrating improved efficiency for filamentation
problems.  Gutnic et al.\ \cite{GutnicBesse2004} used adaptive
phase-space grids with wavelet-based refinement indicators.  The present
work differs in two respects: (i)~we use a DG representation rather than
spline interpolation, preserving the compact two-cell stencil and exact
mass conservation of the SLDG method; and (ii)~we build on PETSc's
p8est-backed AMR infrastructure, which provides scalable parallel mesh
management and enables reuse of the same velocity mesh for both the
transport solver and the Landau collision operator.

\paragraph{Outline.}
Section~\ref{sec:math} presents the mathematical formulation of the SLDG
scheme, the Strang splitting time integrator, and a notation summary
(Table~\ref{tab:notation}).
Section~\ref{sec:amr} describes the AMR extension, including the pencil
data structure, the hybrid fast/slow sweep, and the weighted accumulation
strategy.
Section~\ref{sec:velocity} discusses the 3V velocity space treatment.
Section~\ref{sec:implementation} covers the PETSc implementation.
Section~\ref{sec:results} presents numerical results on Landau damping,
energy conservation, and AMR refinement behavior.
Section~\ref{sec:conclusions} concludes with a summary and directions for
future work.

\section{Mathematical Formulation}
\label{sec:math}

This section presents the Vlasov--Poisson system, the Strang splitting time
integrator, and the SLDG discretization on uniform grids.  The AMR
extension is deferred to Section~\ref{sec:amr}.

\subsection{Vlasov--Poisson system}
\label{sec:vlasov-poisson}

We consider the Vlasov--Poisson system in one spatial dimension with $d_v$
velocity dimensions ($d_v = 1$ or $d_v = 3$).  The distribution function
$f(t,x,\vec{v})$ evolves according to
\begin{align}
  \partial_t f + v_x \,\partial_x f + E(x,t)\,\partial_{v_x} f &= 0,
  \label{eq:vlasov} \\
  -\partial_{xx} \phi &= \rho - 1, \label{eq:poisson} \\
  E &= -\partial_x \phi, \label{eq:efield} \\
  \rho(x,t) &= \int_{\RR^{d_v}} f(x,\vec{v},t) \, \dd\vec{v},
  \label{eq:rho}
\end{align}
where $\vec{v} = (v_x, v_y, v_z)$ in 3V and $\vec{v} = v_x$ in 1V.  The
spatial domain is periodic, $x \in [0, L]$ with $L = 2\pi/k$, and the
velocity domain is $\vec{v} \in [-R, R]^{d_v}$ with absorbing or periodic
boundary conditions.  The background ion charge density is normalized to
unity.

For the Landau damping benchmark, the initial condition is a perturbed
Maxwellian:
\begin{equation}\label{eq:landau-ic}
  f(0,x,\vec{v}) = \frac{1}{(2\pi)^{d_v/2}}
  \exp\!\Bigl(-\tfrac{|\vec{v}|^2}{2}\Bigr)
  \bigl(1 + \alpha \cos(kx)\bigr),
\end{equation}
with perturbation amplitude $\alpha$ and wave number~$k$.  The standard
linear Landau damping test uses $\alpha = 0.01$, $k = 0.5$, and the
theoretical damping rate (the decay rate of the $E_{\max}$ envelope) is
$\gamma = -0.1533$ \cite{LandauDamping1946}.

\subsection{Strang splitting}
\label{sec:strang}

Following Cheng and Knorr \cite{ChengKnorr1976}, we decompose the
Vlasov--Poisson system into a sequence of one-dimensional advections using
Strang splitting \cite{Strang1968}.  Each time step of size $\Delta t$
proceeds as:
\begin{enumerate}
\item \textbf{Half x-advection.}  Solve
  $\partial_t f + v_x \,\partial_x f = 0$ for $\Delta t/2$.
  Each velocity DOF advects independently in~$x$ with speed~$v_x$.
\item \textbf{Field solve.}  Compute the charge density~$\rho$ by
  integrating $f$ over velocity space, solve the Poisson equation
  \eqref{eq:poisson} for~$\phi$, and obtain the electric field
  $E = -\partial_x \phi$.
\item \textbf{Full v-advection.}  Solve
  $\partial_t f + E(x)\,\partial_{v_x} f = 0$ for $\Delta t$.
  Each spatial DOF advects independently in~$v_x$ with speed~$E(x)$.
  In 3V, the $v_y$ and $v_z$ sweeps have zero shift for the
  electrostatic problem and are skipped.
\item \textbf{Half x-advection.}  Repeat step~1 for $\Delta t/2$.
\end{enumerate}
The Strang splitting is second-order accurate in time.  Each advection
substep is a one-dimensional translation that is discretized by the SLDG
method described next.

\paragraph{Time step selection.}
Although the SLDG method is unconditionally stable (no CFL restriction),
the time step $\Delta t$ is limited by three practical considerations.
First, the Strang splitting introduces a splitting error of
$O(\Delta t^2)$, so accuracy requires $\Delta t$ small enough that the
splitting error does not dominate the spatial discretization error.
Second, the ghost stencil width for the x-advection grows with $\Delta t$:
$\mathrm{sw} = \lfloor v_{\max} \Delta t / (2h_x) \rfloor + 2$ (the
factor of~$2$ arises from the half-step $\Delta t/2$ in the Strang
splitting), increasing
communication volume.  Third, for the Landau damping benchmark the
electric field oscillates with period $T \approx 2\pi / \omega_r \approx
4.6$, so $\Delta t$ must resolve these oscillations.  In all experiments
we use $\Delta t = 0.1$, which balances these constraints: the splitting
error is $O(10^{-2})$, the ghost width is $\mathrm{sw} \le 5$ cells for
$R = 6$ and $N_x = 64$, and each oscillation period is resolved by
${\sim}46$ time steps.

\subsection{The SLDG method on uniform grids}
\label{sec:sldg-uniform}

We follow the formulation of Crouseilles, Mehrenberger, and Vecil
\cite{CrouseillesMehrenbergerVecil2011} as implemented by Einkemmer
\cite{Einkemmer2016}.  Consider a one-dimensional advection
$\partial_t u + a\,\partial_\xi u = 0$ on a uniform grid with $N$~cells
of width~$h$ and cell boundaries $\xi_{i-1/2}$, $\xi_{i+1/2}$.  The
distribution is represented by a piecewise polynomial of degree~$p$ on
each cell, with $o = p+1$ degrees of freedom per cell stored at the
Gauss--Lobatto--Legendre (GLL) nodes $\{\hat\xi_j\}_{j=0}^p$ on the
reference element $[-1,1]$.  The corresponding Lagrange basis functions
are $\ell_j(\hat\xi)$.

After a time step $\Delta t$, the exact solution translates the
distribution by $a\Delta t$.  Decomposing the shift into an integer part
and a fractional part,
\begin{equation}\label{eq:shift-decomp}
  \frac{a\Delta t}{h} = n_s + \alpha, \qquad
  n_s = \Bigl\lfloor \frac{a\Delta t}{h} \Bigr\rfloor, \quad
  \alpha \in [0,1),
\end{equation}
the translated distribution in cell~$i$ draws data from cells $i^* = i -
n_s$ (the ``same cell'') and $i^* - 1$ (the ``neighbor cell'').  The $L^2$
projection onto the DG basis yields the update
\begin{equation}\label{eq:sldg-update}
  \tilde{u}^{n+1}_{ij} = \sum_{j'} A^\alpha_{jj'}\,
  \tilde{u}^n_{i^*,j'} + \sum_{j'} B^\alpha_{jj'}\,
  \tilde{u}^n_{i^*-1,j'},
\end{equation}
where the overlap matrices $A^\alpha$ and $B^\alpha$ are given by
\begin{align}
  A^\alpha_{jj'} &= \frac{1}{\omega_j} \sum_{r=0}^{p}
  \omega_r\, \ell_{j'}\!\bigl(\alpha + (1-\alpha)\hat\xi_r\bigr)\,
  \ell_j\!\bigl((1-\alpha)\hat\xi_r\bigr),
  \label{eq:A-matrix} \\
  B^\alpha_{jj'} &= \frac{1}{\omega_j} \sum_{r=0}^{p}
  \omega_r\, \ell_{j'}\!\bigl(\alpha\hat\xi_r\bigr)\,
  \ell_j\!\bigl(\alpha(\hat\xi_r - 1) + 1\bigr),
  \label{eq:B-matrix}
\end{align}
with $\omega_j$ the GLL quadrature weights and $\hat\xi_r$ the Gauss
quadrature points on $[0,1]$ (distinct from the GLL nodes
$\{\hat\xi_j\}$ used as DOF locations).  For $p \le 1$ (piecewise constant or
linear), the GLL mass matrix coincides with the lumped (diagonal) mass
matrix and the division by $\omega_j$ suffices.  For $p \ge 2$, the
mass matrix $M_{jj'} = \int_{-1}^{1}
\ell_j(\hat\xi)\,\ell_{j'}(\hat\xi)\,\dd\hat\xi$ is not diagonal and
must be inverted:
\begin{equation}\label{eq:minv}
  A^\alpha \leftarrow M^{-1} A^\alpha, \qquad
  B^\alpha \leftarrow M^{-1} B^\alpha.
\end{equation}

\paragraph{Conservation.}
The SLDG method is conservative by construction: the overlap integrals
satisfy the partition-of-unity property
\begin{equation}\label{eq:partition-unity}
  \sum_i \omega_i \bigl(A^\alpha_{ij} + B^\alpha_{ij}\bigr) = \omega_j
  \qquad \text{for all } j,
\end{equation}
where $\omega_i$ are the GLL quadrature weights.  This ensures that the
total mass $\sum_{i,j} \omega_j\, \tilde{u}_{ij}$ is preserved exactly
(to machine precision) at each advection step.

\paragraph{Key property: two-cell stencil.}
A crucial advantage of the SLDG method is that the update
\eqref{eq:sldg-update} requires data from exactly two source cells ($i^*$
and $i^*-1$), regardless of the polynomial degree~$p$.  This is in
contrast to polynomial interpolation methods, where the stencil width
grows with the order of accuracy.  The compact stencil minimizes
communication overhead in parallel implementations and, as we show in
Section~\ref{sec:amr}, enables a natural extension to AMR grids where the
two source cells may differ in size
(Fig.~\ref{fig:sldg-amr-geometry}).

\subsection{Precomputation and the uniform-grid assumption}
\label{sec:precompute}

On a uniform grid, the fractional shift $\alpha$ depends only on the
advection speed and the cell width, not on the cell index.  For the
x-advection step, each velocity DOF has a different speed $v_x$, yielding
a different $\alpha$ value, but the matrices $A^\alpha$ and $B^\alpha$ can
be precomputed for each velocity DOF before the time loop begins.  For the
v-advection step, the electric field $E(x)$ varies with position, so the
matrices must be recomputed at each time step, but on a uniform grid all
cells at the same refinement level share the same cell width~$h$, and the
precomputation cost is $O(p^2)$ per distinct shift value.

This precomputation strategy breaks down on an AMR grid, where cells at
different refinement levels have different widths $h_\ell = h_0 / 2^\ell$.
The fractional shift $\alpha_\ell = (a\Delta t / h_\ell) - \lfloor
a\Delta t / h_\ell \rfloor$ differs by level, and at refinement
boundaries the source and destination cells may have different sizes,
requiring a generalized overlap integral.  This extension is the subject
of Section~\ref{sec:amr}.

\subsection{Notation summary}
\label{sec:notation}

Table~\ref{tab:notation} collects the principal symbols used throughout
the paper.

\begin{table}[htbp]
\centering
\caption{Summary of notation.}
\label{tab:notation}
\smallskip
\begin{tabular}{cl}
\toprule
Symbol & Meaning \\
\midrule
$f(t,x,\vec{v})$ & Distribution function \\
$\alpha$ & Fractional shift in SLDG update \\
$p$ & Polynomial degree of DG basis \\
$p_x$ & Polynomial degree in x-space \\
$h$, $h_\ell$ & Cell width (at refinement level $\ell$) \\
$\ell$, $L$ & Refinement level; maximum AMR level \\
$N_b$ & Number of base cells per velocity dimension \\
$N_x$ & Number of x-space cells \\
$R$ & Velocity domain radius ($\vec{v} \in [-R,R]^3$) \\
$k$ & Wave number of initial perturbation \\
$\gamma$ & Landau damping rate \\
$A^\alpha$, $B^\alpha$ & Same-cell and neighbor-cell overlap matrices \\
$M^{(s,c_s)}$ & Generalized overlap matrix (slow path) \\
$\omega_j$ & GLL quadrature weight at node $j$ \\
$w_c$ & Transverse pencil weight for cell $c$ \\
$n_c$ & Number of pencils containing cell $c$ \\
\bottomrule
\end{tabular}
\end{table}

\section{AMR Extension}
\label{sec:amr}

This section describes the extension of the SLDG method to adaptively
refined velocity grids.  The mesh is generated by PETSc's
\texttt{DMPlexLandauCreateVelocitySpace} using p4est (2D) or p8est (3D)
\cite{BursteddeWilcoxGhattas2011} as the forest-of-octrees back end.  The
AMR strategy is a heuristic inherited from PETSc's Landau collision
operator \cite{AdamsLandau2025}: cells are refined near the origin of
velocity space to resolve the shoulders of the Maxwellian distribution
(where the gradient $|\nabla_v f|$ is largest, at $|v| \sim v_{\rm th}$)
on the coarse grids used by the collision operator, while cells
near the domain boundary are left coarse because the distribution is
exponentially small there.

\subsection{Cell geometry and refinement levels}
\label{sec:cell-geometry}

After mesh creation, each cell~$c$ is assigned a refinement level
$\ell_c \in \{0, 1, \ldots, L\}$ and its physical extent is recorded:
the lower bound $v^{\mathrm{lo}}_{c,d}$ and width $h_{c,d}$ in each
coordinate direction $d \in \{0,1,2\}$.  On a mesh with $L$~levels of
refinement, the cell width at level~$\ell$ in direction~$d$ is
\begin{equation}\label{eq:h-level}
  h_{\ell,d} = \frac{h_{0,d}}{2^\ell},
\end{equation}
where $h_{0,d} = 2R / N_{0,d}$ is the coarsest cell width, $R$ is the
velocity domain radius, and $N_{0,d}$ is the number of coarse cells in
direction~$d$.

\subsection{Pencil data structure}
\label{sec:pencils}

The SLDG sweep is inherently one-dimensional: for each sweep direction~$d$,
the update \eqref{eq:sldg-update} operates along lines of cells aligned
with the $d$-axis.  On a uniform grid, these lines are trivially defined
by the Cartesian structure.  On an AMR grid, the lines must be extracted
from the unstructured mesh topology.

We define a \emph{pencil} as a maximal sequence of cells that share the
same transverse extent and tile the sweep axis without gaps.  For a 3D
mesh with sweep direction~$d$, the two transverse dimensions define a 2D
grid of intervals.  Each pair of transverse intervals $(I_{t_1}, I_{t_2})$
defines one pencil: a cell belongs to this pencil if its transverse extent
covers both intervals.

The pencil extraction algorithm proceeds as follows:
\begin{enumerate}
\item Collect all unique edge coordinates in each transverse dimension
  from the cell bounding boxes.
\item Form the Cartesian product of the resulting 1D interval sets to
  obtain a 2D grid of transverse intervals.
\item For each transverse interval pair, collect all cells whose
  transverse extent covers the interval, sort them by sweep-direction
  coordinate, and validate gap-free coverage.
\end{enumerate}
The pencils are stored in compressed sparse row (CSR) format: an offset
array indexes into arrays of cell identifiers, cell widths, and cell lower
bounds along the sweep direction.

\paragraph{Coarse cells in multiple pencils.}
On a uniform mesh, each cell appears in exactly one pencil per sweep
direction.  On an AMR mesh, a coarse cell at level~$\ell$ may span
multiple fine-level transverse intervals and therefore appear in up to
$4^{L-\ell}$ pencils (in 3D).  Each pencil entry records a transverse
weight $w_c = 1/n_c$, where $n_c$ is the number of pencils containing
cell~$c$ in this sweep direction.  This is illustrated in
Fig.~\ref{fig:sldg-amr-2v}, where the coarse cell~$L$ spans both fine
transverse intervals and appears in $n_c = 2$ pencils.  The weighted
accumulation strategy for the write-back is described in
Section~\ref{sec:writeback}.

\subsection{Conforming and nonconforming cells}
\label{sec:classify}

Within each pencil, a cell is classified as \emph{conforming} if all
neighbors within a distance of $\pm 2$ cells in the pencil are at the
same refinement level, and \emph{nonconforming} otherwise.  (The radius
of~2 ensures that both source cells of the SLDG stencil, and their own
source cells, are at the same level.)  On a uniform
mesh ($L = 0$), all cells are conforming.  On an AMR mesh, nonconforming
cells occur at refinement boundaries where the cell width changes.

\begin{figure}[htbp]
\centering
\includegraphics[width=0.82\textwidth]{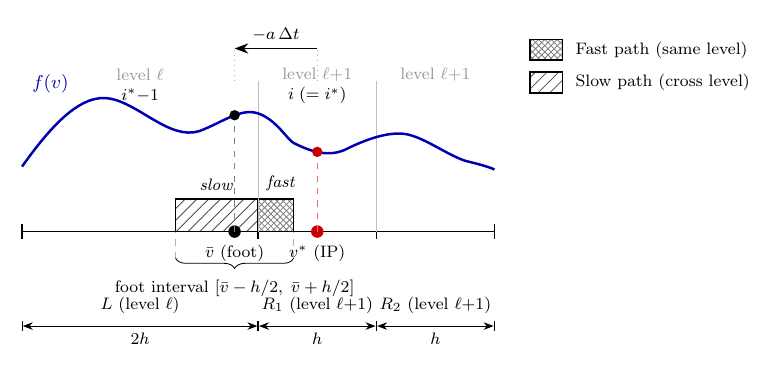}
\caption{SLDG sweep geometry on an AMR mesh (1D cross-section).
  The mesh has a coarse cell~$L$ (width $2h$, level~$\ell$) and two fine
  cells $R_1$, $R_2$ (width $h$, level~$\ell{+}1$).  The red dot (right) marks
  the target IP~$v^*$; the black dot marks the foot point
  $\bar{v} = v^* - a\,\Delta t$.  The brace shows the foot interval.
  Shaded rectangles indicate the fast-path (crosshatch,
  same-level source) and slow-path (NE hatching, cross-level source)
  contributions to the $L^2$ projection.}
\label{fig:sldg-amr-geometry}
\end{figure}

This classification determines which sweep kernel is used for each cell:
\begin{itemize}
\item \textbf{Fast path} (conforming cells): precomputed per-level
  $A^\alpha$ and $B^\alpha$ matrices, identical to the uniform-grid
  method.
\item \textbf{Slow path} (nonconforming cells): generalized overlap
  integrals computed on the fly, handling arbitrary source/destination
  cell sizes.
\end{itemize}
Figure~\ref{fig:sldg-amr-geometry} illustrates this classification for a
fine-cell IP whose foot interval straddles the coarse--fine boundary;
Fig.~\ref{fig:sldg-amr-2v} extends this to 2D, showing both sweep
directions simultaneously.

\paragraph{Sweep geometry.}
Consider a target integration point (IP) $v^*$, one of the $p{+}1$ GLL
quadrature nodes inside a fine cell~$R_1$ at level~$\ell{+}1$
(Fig.~\ref{fig:sldg-amr-geometry}).  The updated value
$f^{\mathrm{out}}_{R_1}(v^*)$ is obtained by tracing the characteristic
backward by one step to the \emph{foot point}
$\bar{v} = v^* - a\,\Delta t$.  The \emph{foot interval}
$[\bar{v} - h_{R_1}/2,\; \bar{v} + h_{R_1}/2]$ is the physical extent of
cell~$R_1$ shifted upstream by $a\,\Delta t$; it is the support of the
translated polynomial that must be projected back onto the DG basis.

When the foot interval straddles a refinement boundary, it overlaps source
cells at different levels.  The $L^2$
projection~\eqref{eq:sldg-update} is split into two contributions.  The
portion of the foot interval that lies inside a same-level source cell
(the \emph{fast-path chunk}, crosshatch in
Fig.~\ref{fig:sldg-amr-geometry}) is handled by the precomputed matrix
$A^{\alpha_{\ell+1}}$~\eqref{eq:fast-update}.  The portion that overlaps
a coarser cell~$L$ at level~$\ell$ (the \emph{slow-path chunk}, NE
hatching) requires the generalized overlap integral
$M^{(R_1,L)}$~\eqref{eq:slow-overlap}.  Together the two chunks cover the
entire foot interval, ensuring that no mass is lost or gained (partition of
unity, eq.~\eqref{eq:partition-unity}).

\subsection{Per-level precomputation (fast path)}
\label{sec:fast-path}

For each refinement level $\ell$ and a given advection speed~$a$, the
shift in level-$\ell$ cells is
\begin{equation}\label{eq:shift-level}
  \frac{a\,\Delta t}{h_{\ell,d}} = n_{s,\ell} + \alpha_\ell, \qquad
  n_{s,\ell} = \Bigl\lfloor \frac{a\,\Delta t}{h_{\ell,d}} \Bigr\rfloor,
  \quad \alpha_\ell \in [0,1).
\end{equation}
The overlap matrices $A^{\alpha_\ell}$ and $B^{\alpha_\ell}$ are computed
from \eqref{eq:A-matrix}--\eqref{eq:B-matrix} exactly as on a uniform
grid.  For a conforming cell at position~$s$ in the pencil, the source
cells are at positions $s - n_{s,\ell}$ (same cell) and
$s - n_{s,\ell} - 1$ (neighbor), and the update is
\begin{equation}\label{eq:fast-update}
  f^{\mathrm{out}}_{s,i} = \sum_j A^{\alpha_\ell}_{ij}\,
  f^{\mathrm{in}}_{s-n_{s,\ell},j} + \sum_j B^{\alpha_\ell}_{ij}\,
  f^{\mathrm{in}}_{s-n_{s,\ell}-1,j}.
\end{equation}
This is $O(p^2)$ work per cell, identical to the uniform-grid cost.
The fast-path region corresponds to the horizontally hatched chunk in
Fig.~\ref{fig:sldg-amr-geometry}.

\subsection{Generalized overlap integral (slow path)}
\label{sec:slow-path}

At refinement boundaries, the source and destination cells may have
different widths.  The foot of the characteristic from destination cell~$s$
spans the physical interval $[v^{\mathrm{lo}}_s - a\Delta t,\;
v^{\mathrm{lo}}_s + h_s - a\Delta t]$, which may overlap with multiple
source cells of varying sizes.  The generalized overlap integral computes
the contribution from each overlapping source cell~$c_s$:
\begin{equation}\label{eq:slow-overlap}
  M_{ij}^{(s,c_s)} = \frac{1}{h_s/2} \int_{v_L}^{v_R}
  \ell_i\!\Bigl(\frac{2(v + a\Delta t - v^{\mathrm{lo}}_s)}{h_s} - 1\Bigr)\,
  \ell_j\!\Bigl(\frac{2(v - v^{\mathrm{lo}}_{c_s})}{h_{c_s}} - 1\Bigr)
  \,\dd v,
\end{equation}
where $[v_L, v_R]$ is the intersection of the foot interval with source
cell~$c_s$, and the arguments to $\ell_i$ and $\ell_j$ map the physical
coordinate to the reference element $[-1,1]$ of the destination and source
cells, respectively.  The integral is evaluated by Gauss quadrature with
$2p+2$ points, mapped to the overlap interval $[v_L, v_R]$.

After applying the mass matrix inverse \eqref{eq:minv}, the contribution
from source cell~$c_s$ to destination cell~$s$ is accumulated:
\begin{equation}\label{eq:slow-update}
  f^{\mathrm{out}}_{s,i} \mathrel{+}= \sum_j
  (M^{-1} M^{(s,c_s)})_{ij}\, f^{\mathrm{in}}_{c_s,j}.
\end{equation}
The slow path is $O(p^2 \cdot n_{\mathrm{src}})$ per cell, where
$n_{\mathrm{src}}$ is the number of overlapping source cells (typically 2
at a refinement boundary).  The slow-path region corresponds to the
NE-hatched chunk in Fig.~\ref{fig:sldg-amr-geometry}.  On a uniform mesh,
the slow path reduces to the fast path with $n_{\mathrm{src}} = 2$.

\paragraph{Absorbing boundary conditions.}
At the velocity domain boundary, source cells outside the domain contribute
zero to the overlap integral, implementing absorbing boundary conditions.

\paragraph{Periodic boundary conditions.}
Periodic boundary conditions in velocity space are also supported as an
alternative.  When the foot of a characteristic falls outside the domain,
the source cell index is wrapped modulo the pencil length before the
overlap integral is evaluated, so the update formula \eqref{eq:slow-update}
is otherwise unchanged.
This maintains conservation of mass.

\subsection{Weighted accumulation for AMR write-back}
\label{sec:writeback}

\begin{figure}[htbp]
\centering
\includegraphics[width=0.82\textwidth]{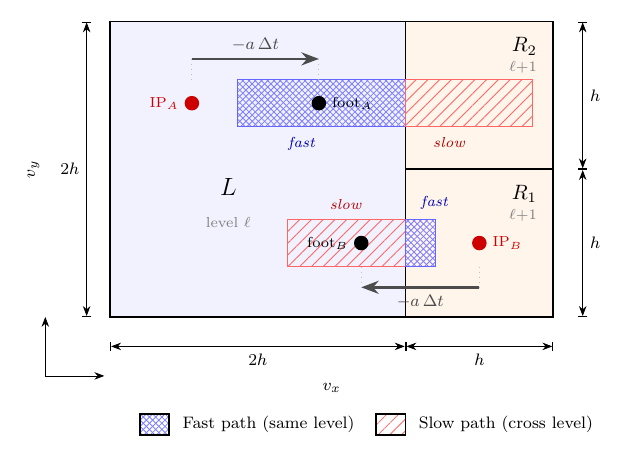}
\caption{2D SLDG sweep geometry on an AMR mesh.  Coarse cell~$L$
  (width~$2h$) and fine cells~$R_1$, $R_2$ (width~$h$).  Upper band
  (Sweep~A): coarse-cell IP sweeps rightward into~$R_2$.  Lower band
  (Sweep~B): fine-cell IP sweeps leftward into~$L$.  Crosshatch
  = fast path (same level); NE hatching = slow path (cross
  level).  The vertical extent of the bands is schematic; the overlap
  integrals are 1D (see Fig.~\ref{fig:sldg-amr-geometry}).}
\label{fig:sldg-amr-2v}
\end{figure}

On an AMR mesh, a coarse cell~$c$ at level~$\ell < L$ appears in multiple
pencils with transverse weight $w_c = 1/n_c$.  The write-back from pencil
results to the global solution uses weighted accumulation:
\begin{equation}\label{eq:weighted-writeback}
  f^{\mathrm{out}}_c = f^{\mathrm{in}}_c + \sum_{p \ni c}
  w_c \bigl(f^{\mathrm{pencil},p}_c - f^{\mathrm{in}}_c\bigr),
\end{equation}
where the sum is over all pencils~$p$ containing cell~$c$, and
$f^{\mathrm{pencil},p}_c$ is the result of the SLDG sweep in pencil~$p$.
For a cell appearing in $n_c$ pencils, the weights sum to unity:
$\sum_{p \ni c} w_c = n_c \cdot (1/n_c) = 1$, so the final result is the
average of the pencil results.  On a uniform mesh where $n_c = 1$, this
reduces to a simple copy: $f^{\mathrm{out}}_c = f^{\mathrm{pencil}}_c$.
Figure~\ref{fig:sldg-amr-2v} illustrates the pencil weights for a 2D AMR
mesh.  The figure shows two representative pencil sweeps in the $v_x$
direction on a mesh with one coarse cell~$L$ (level~$\ell$, width~$2h$)
and two fine cells~$R_1$, $R_2$ (level~$\ell{+}1$, width~$h$).
In \textbf{Sweep~A}, the target IP is a GLL node in the coarse cell~$L$
and the characteristic traces rightward into fine cell~$R_2$; the foot
interval straddles the coarse--fine boundary, producing a fast-path chunk
in~$L$ (precomputed $A^{\alpha_\ell}$) and a slow-path chunk in~$R_2$
(generalized overlap integral $M^{(s,c_s)}$).
In \textbf{Sweep~B}, the target IP is in fine cell~$R_1$ and the
characteristic traces leftward into~$L$; the roles are reversed, with the
fast-path chunk in~$R_1$ and the slow-path chunk in~$L$.
In 2D, cell~$L$ spans both fine transverse intervals and therefore appears
in $n_c = 2$ pencils with weight $w_c = 1/2$ each, while $R_1$ and $R_2$
each appear in one pencil ($n_c = 1$, $w_c = 1$).  The weighted
accumulation~\eqref{eq:weighted-writeback} averages the pencil results
for~$L$, ensuring that the total weight sums to unity.

\subsection{Partition of unity on AMR grids}
\label{sec:amr-conservation}

The fast-path update \eqref{eq:fast-update} inherits the partition-of-unity
property \eqref{eq:partition-unity} from the uniform-grid method, since it
uses the same $A^\alpha$ and $B^\alpha$ matrices.  The slow-path update
\eqref{eq:slow-update} also preserves mass: the generalized overlap
integral \eqref{eq:slow-overlap} accounts for the full foot interval, and
the contributions from all source cells sum to the complete $L^2$
projection.  The weighted accumulation \eqref{eq:weighted-writeback}
preserves mass because the weights sum to unity for each cell.  We verify
this numerically in Section~\ref{sec:results}.
Algorithm~\ref{alg:sldg-sweep} summarizes the hybrid SLDG sweep for one
velocity dimension on an AMR mesh.

\begin{algorithm}[htbp]
\caption{Hybrid SLDG sweep along one velocity dimension on an AMR mesh.}
\label{alg:sldg-sweep}
\begin{algorithmic}[1]
\Require Distribution $f^{\mathrm{in}}$, advection speed $a$, time step
  $\Delta t$, pencil set $\mathcal{P}_d$ for sweep direction $d$
\Ensure Updated distribution $f^{\mathrm{out}}$
\State Precompute $A^{\alpha_\ell}$, $B^{\alpha_\ell}$ for each
  refinement level $\ell$ \Comment{Eq.~\eqref{eq:shift-level}}
\State Initialize $f^{\mathrm{out}} \gets f^{\mathrm{in}}$
\ForAll{pencils $p \in \mathcal{P}_d$}
  \State \textsc{Pack}: extract 1D pencil data from $f^{\mathrm{in}}$
    using tensor-product inverse map
  \ForAll{cells $s$ in pencil $p$}
    \If{cell $s$ is \emph{conforming} (same-level neighbors within $\pm 2$)}
      \State \textbf{Fast path}: apply precomputed matrices
        \Comment{Eq.~\eqref{eq:fast-update}}
      \State $f^{\mathrm{pencil}}_s \gets A^{\alpha_\ell} f^{\mathrm{in}}_{s-n_s}
        + B^{\alpha_\ell} f^{\mathrm{in}}_{s-n_s-1}$
    \Else
      \State \textbf{Slow path}: compute generalized overlap integrals
        \Comment{Eq.~\eqref{eq:slow-overlap}}
      \State $f^{\mathrm{pencil}}_s \gets \sum_{c_s} M^{-1} M^{(s,c_s)}
        f^{\mathrm{in}}_{c_s}$
    \EndIf
  \EndFor
  \State \textsc{Scatter}: weighted accumulation to $f^{\mathrm{out}}$
    with weight $w_c = 1/n_c$ \Comment{Eq.~\eqref{eq:weighted-writeback}}
\EndFor
\end{algorithmic}
\end{algorithm}

\subsection{AMR velocity mesh geometry}
\label{sec:amr-mesh-geometry}

\begin{figure}[htbp]
\centering
\begin{subfigure}[t]{0.41\textwidth}
  \centering
  \includegraphics[width=\textwidth]{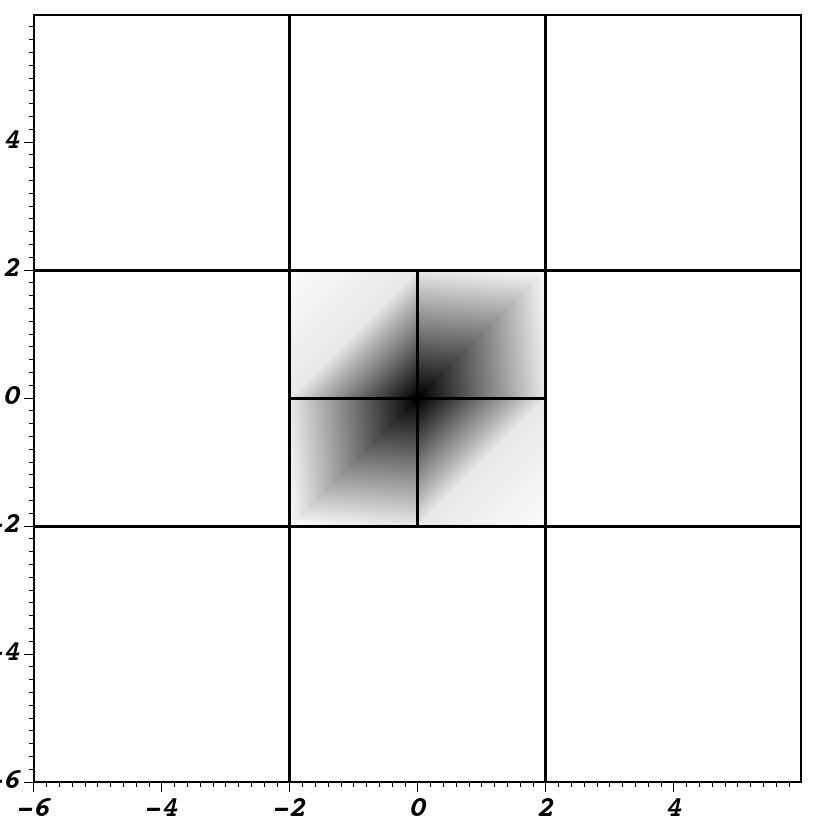}
  \caption{$3^3$+AMR1 (34 cells)}
  \label{fig:amr-mesh-333}
\end{subfigure}
\hfill
\begin{subfigure}[t]{0.51\textwidth}
  \centering
  \includegraphics[width=\textwidth]{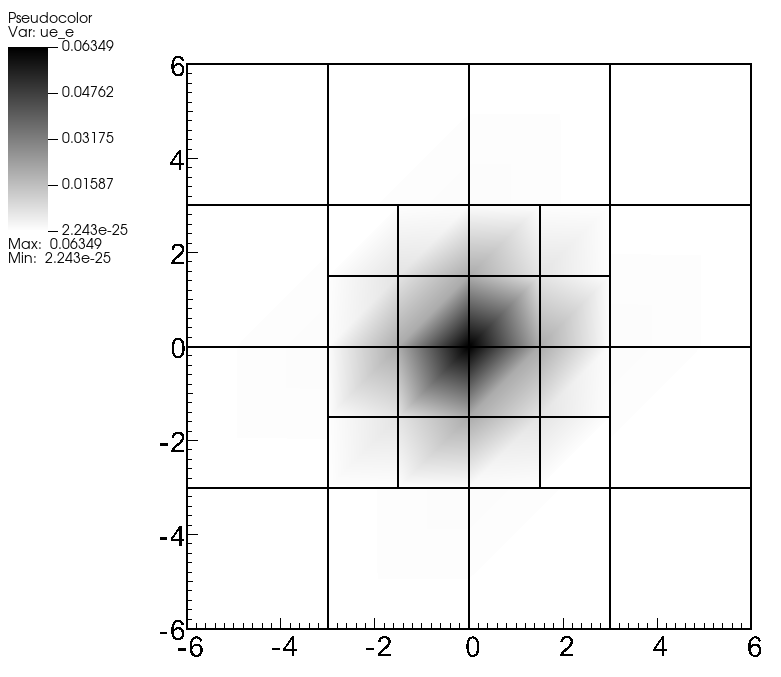}
  \caption{$4^3$+AMR1 (120 cells)}
  \label{fig:amr-mesh-444}
\end{subfigure}
\caption{2D cross-sections ($v_x$--$v_y$ plane at $v_z = 0$) of the
  electron distribution $u_e$ on two AMR velocity meshes ($R=6$).
  (a)~The $3^3$ base mesh with one AMR level; the four inner base cells
  are each refined once, concentrating resolution near the Maxwellian
  peak.  (b)~The $4^3$ base mesh with one AMR level; the finer base
  grid resolves the Maxwellian shoulders more accurately. (Visit uses linear interpolation from vertices)}
\label{fig:amr-meshes}
\end{figure}

The velocity mesh is created by
\texttt{DMPlexLandauCreateVelocitySpace} with an $N_b^3$ base grid and
up to $L$ levels of p8est AMR refinement about the origin.  The AMR
strategy refines the $2^3 = 8$ cells nearest the origin, concentrating
resolution near the Maxwellian peak while keeping coarse cells where the
distribution is exponentially small.
Figure~\ref{fig:amr-meshes} shows representative 2D cross-sections for
the $3^3$ and one level of mesh adaptivity (AMR1) and $4^3$+AMR1 configurations; cell counts for all
configurations are reported in Table~\ref{tab:convergence}.

\section{3V Velocity Space Treatment}
\label{sec:velocity}

The extension from 1V to 3V velocity space requires two ingredients:
tensor-product DG elements on hexahedral cells, and a dimensional splitting
strategy that reduces the 3D advection to a sequence of 1D pencil sweeps.

\subsection{Tensor-product DG elements}
\label{sec:tensor-product}

The velocity mesh consists of hexahedral cells, each equipped with a
tensor-product DG basis of degree~$p$ in each direction.  The 3D basis
functions are products of 1D Lagrange polynomials at GLL nodes:
\begin{equation}\label{eq:tp-basis}
  \phi_{(i,j,k)}(\hat\xi, \hat\eta, \hat\zeta) =
  \ell_i(\hat\xi)\,\ell_j(\hat\eta)\,\ell_k(\hat\zeta),
\end{equation}
where $(\hat\xi, \hat\eta, \hat\zeta) \in [-1,1]^3$ are reference
coordinates and $i,j,k \in \{0,\ldots,p\}$.  Each cell has
$N_{\rm dof} = (p+1)^3$ degrees of freedom.  The tensor-product structure is
essential: it allows the 3D advection to be decomposed into independent 1D
sweeps along each coordinate axis.

\subsection{Tensor-product permutation}
\label{sec:tp-permutation}

PETSc's finite element infrastructure (PetscFE) assigns DOF indices to
each cell in an order that may not align with the tensor-product structure.
We discover the tensor-product permutation by tabulating the basis
functions at the GLL grid points: for each DOF~$b$, we find the unique
grid point $(\hat\xi_i, \hat\eta_j, \hat\zeta_k)$ where
$\phi_b \approx 1$, yielding the 1D indices $(i,j,k)$.

For each sweep direction $d \in \{0,1,2\}$, we define:
\begin{itemize}
\item the \emph{sweep index}: the 1D index along direction~$d$,
\item the \emph{transverse indices}: the two 1D indices in the
  perpendicular directions.
\end{itemize}
An inverse map is precomputed so that, given a sweep direction and a pair
of transverse indices $(t_1, t_2)$, the 1D pencil data can be packed and
unpacked in $O(p+1)$ operations per cell rather than scanning all
$N_{\rm dof} = (p+1)^3$ DOFs.

\subsection{Dimensional splitting in 3V}
\label{sec:dim-splitting}

The v-advection step solves $\partial_t f + E(x)\,\partial_{v_x} f = 0$
for the full time step $\Delta t$.  In the electrostatic 1X+3V problem,
only the $v_x$ component is driven by the electric field; the $v_y$ and
$v_z$ sweeps have zero shift and are skipped.  For the general
electromagnetic case, all three sweeps would be active with shifts
determined by the Lorentz force components.

For each x-space DOF $(c_x, b_x)$, the electric field $E(x)$ is
interpolated at the GLL node position to obtain the advection speed.  The
3D v-advection then proceeds as a sequence of 1D pencil sweeps:
\begin{enumerate}
\item For each pair of transverse indices $(t_1, t_2)$ and each pencil~$p$
  in the sweep direction:
  \begin{enumerate}
  \item \textbf{Pack}: extract the 1D pencil data
    $f_{\mathrm{pencil}}[s \cdot (p+1) + b_{\mathrm{sweep}}]$ from the
    global array using the tensor-product inverse map.
  \item \textbf{Sweep}: apply the hybrid fast/slow SLDG update
    (Section~\ref{sec:amr}) along the pencil.
  \item \textbf{Scatter}: write the results back to the global array
    using weighted accumulation \eqref{eq:weighted-writeback}.
  \end{enumerate}
\item Repeat for each active sweep direction.
\end{enumerate}
The pack and scatter steps use the precomputed inverse map from
Section~\ref{sec:tp-permutation}, making them $O(p+1)$ per cell rather
than $O((p+1)^3)$.

\subsection{X-space advection}
\label{sec:x-advection}

The x-advection step solves $\partial_t f + v_x\,\partial_x f = 0$ for
$\Delta t/2$.  The spatial domain is discretized on a 1D periodic grid
(DMDA) with $N_x$ cells and DG degree $p_x$.  Each velocity DOF has a
different advection speed~$v_x$, so the SLDG overlap matrices $A^\alpha$
and $B^\alpha$ are precomputed for each velocity DOF before the time loop.

The x-advection requires ghost data from neighboring MPI processes.  To
minimize communication, all velocity DOFs at each x-cell are packed into a
single batched vector with $N_v \cdot (p_x+1)$ components per cell, and a
single \texttt{DMGlobalToLocal} call exchanges all ghost data
simultaneously.  The SLDG update is then applied independently for each
velocity DOF using the precomputed matrices.

\section{PETSc Implementation}
\label{sec:implementation}

The solver is implemented as a PETSc example using the DMPlex and DMDA
infrastructure.  The key design decisions are summarized here.

\paragraph{Velocity mesh.}
The velocity-space mesh is created by
\texttt{DMPlexLandauCreateVelocitySpace}, which generates a 3D hexahedral
mesh with optional p8est-backed AMR.  The SLDG solver clones this mesh and
attaches a DG finite element space with the \texttt{v\_petscspace\_degree}
option.  This reuse of the Landau mesh infrastructure ensures that the
same AMR refinement strategy (refinement near the origin, coarsening at
the boundary) is available for both the collision operator and the
transport solver.

\paragraph{Spatial mesh.}
The x-space uses a 1D periodic DMDA with ghost width determined by the
maximum CFL number: $\mathrm{sw} = \lfloor v_{\max} \Delta t / (2h_x)
\rfloor + 2$.

\paragraph{Poisson solver.}
The 1D Poisson equation is discretized with high-order ($p \ge 2$) finite
elements on the DMDA grid and solved with the conjugate gradient (CG) method
preconditioned by GAMG, PETSc's built-in smoothed aggregation algebraic
multigrid solver.  The null space arising from periodic boundary conditions
is accommodated with an SVD coarse-grid solver.

\paragraph{Data layout.}
The distribution function is stored as a single PETSc Vec of size
$N_v^{\mathrm{DOF}} \times N_x^{\mathrm{DOF}}$, with velocity-major
ordering: $f[\mathrm{iv} \cdot N_x^{\mathrm{DOF}} + \mathrm{ix}]$.  This
layout provides unit-stride access for the x-advection sweep and strided
access for the v-advection sweep.

\paragraph{Diagnostics.}
At each time step, the solver computes and records the maximum electric
field $E_{\max}$, the zeroth through second velocity moments (mass,
momentum, energy), the electric field energy
$\mathcal{E}_E = \frac{1}{2} h_x \|E\|^2$, and the total conserved
energy $\mathcal{E} = \frac{1}{2} m_2 + \mathcal{E}_E$.  The Landau
damping rate is extracted by fitting $E_{\max}$ on a log scale versus time
to the envelope peaks.

\paragraph{GPU execution via Kokkos.}
The solver uses Kokkos \cite{TrottKokkos2022} for portable GPU execution.
All compute-intensive kernels---the v-advection pencil sweeps (AdvectV),
the x-advection SLDG update (AdvectX), the charge density reduction, and
the moment diagnostics---are implemented as
\texttt{Kokkos::parallel\_for} or \texttt{parallel\_reduce} loops,
executing on the default Kokkos backend (CUDA on NVIDIA GPUs, Serial on
CPU-only builds).  The device execution model and data-transfer
optimizations are described in Section~\ref{sec:gpu-model}.

\section{Numerical Results}
\label{sec:results}

All experiments use the \texttt{ex6} PETSc example compiled with the
optimized build (\texttt{PETSC\_ARCH=arch-perlmutter-opt-gcc-kokkos-cuda-3d})
and run on Perlmutter GPU nodes at NERSC (NVIDIA A100 40\,GB,
4 GPUs per node) using up to 16 MPI processes (one GPU per process).

\subsection{Test problem: linear Landau damping in 1X+3V}
\label{sec:landau-damping}

We verify the solver on the standard linear Landau damping benchmark
\cite{LandauDamping1946} in 1X+3V, following the PETSc example of
Finn \cite{FinnLandauDamping2023}.  The initial distribution function
is the perturbed Maxwellian~\eqref{eq:landau-ic} with perturbation
amplitude $\alpha = 0.01$ and wave number $k = 0.5$ on the periodic
spatial domain $[0, 2\pi/k]$.  The velocity domain is $[-R, R]^3$
with $R = 6$, which captures all but $\sim\!10^{-8}$ of the Maxwellian
mass.  The Vlasov--Poisson system is closed by the 1D Poisson
equation~\eqref{eq:poisson} with charge density
$\rho(x) = \int f \,\dd^3 v$ and periodic boundary conditions.
The electric field $E = -\partial_x \phi$ drives the $v_x$ advection step.

For linear damping ($\alpha \ll 1$), linear theory predicts exponential
decay of the electric field amplitude at the Landau damping rate
$\gamma_{\rm th} = -0.1533$ (for $k = 0.5$, Maxwellian equilibrium
\cite{LandauDamping1946}).  The numerical damping rate is extracted by
fitting $\log E_{\max}(t)$ to the envelope of local maxima (peaks)
before the recurrence onset.

All convergence tests use $N_x = 64$ x-cells, DG degree $p_x = 2$ in
x-space, and time step $\Delta t = 0.1$ (200 steps to $t = 20$).  The
velocity mesh is created by \texttt{DMPlexLandauCreateVelocitySpace}
with an $N_b^3$ base grid and optional p8est AMR with $L$ levels of
refinement about the origin (see Section~\ref{sec:amr-mesh-geometry}).
The number of velocity cells depends on $N_b$ and $L$.  In 3D, the
PETSc AMR strategy refines the $2^3 = 8$ cells nearest the origin at
each level; each refined cell produces $2^3 = 8$ children, so each AMR
level adds $8 \times 8 - 8 = 56$ cells.  For even $N_b$ this gives
$N_{\rm cells} = N_b^3 + 56 L$; for odd $N_b$ the origin does not
coincide with a cell corner, so p8est selects a single cell (e.g., $3^3$+AMR1 yields 34~cells rather than
$27 + 56 = 83$; $3^3$+AMR2 yields 272~cells because the second level
refines a larger neighborhood).  Exact cell counts for all configurations
are reported in Table~\ref{tab:convergence}.  The number of
integration points (IPs) is $N_{\rm IP} = N_{\rm cells} \cdot (p+1)^3$,
where $p$ is the DG degree in velocity space.

\subsection{Convergence}
\label{sec:convergence}

We study convergence of the Landau damping rate along two axes:
\emph{h-refinement} (increasing the number of uniform velocity cells
$N_b^3$ at fixed polynomial degree~$p$) and \emph{p-refinement}
(increasing $p$ at fixed mesh).  We then compare uniform meshes with
AMR meshes at similar integration point (IP) counts to assess AMR
efficiency.  The Landau damping rate is extracted by fitting
$\log E_{\max}$ at the detected peak values (envelope fit) using least
squares; the fit is restricted to the monotonically-decreasing run of
peaks before recurrence onset.

\begin{figure}[htbp]
\centering
\includegraphics[width=\textwidth]{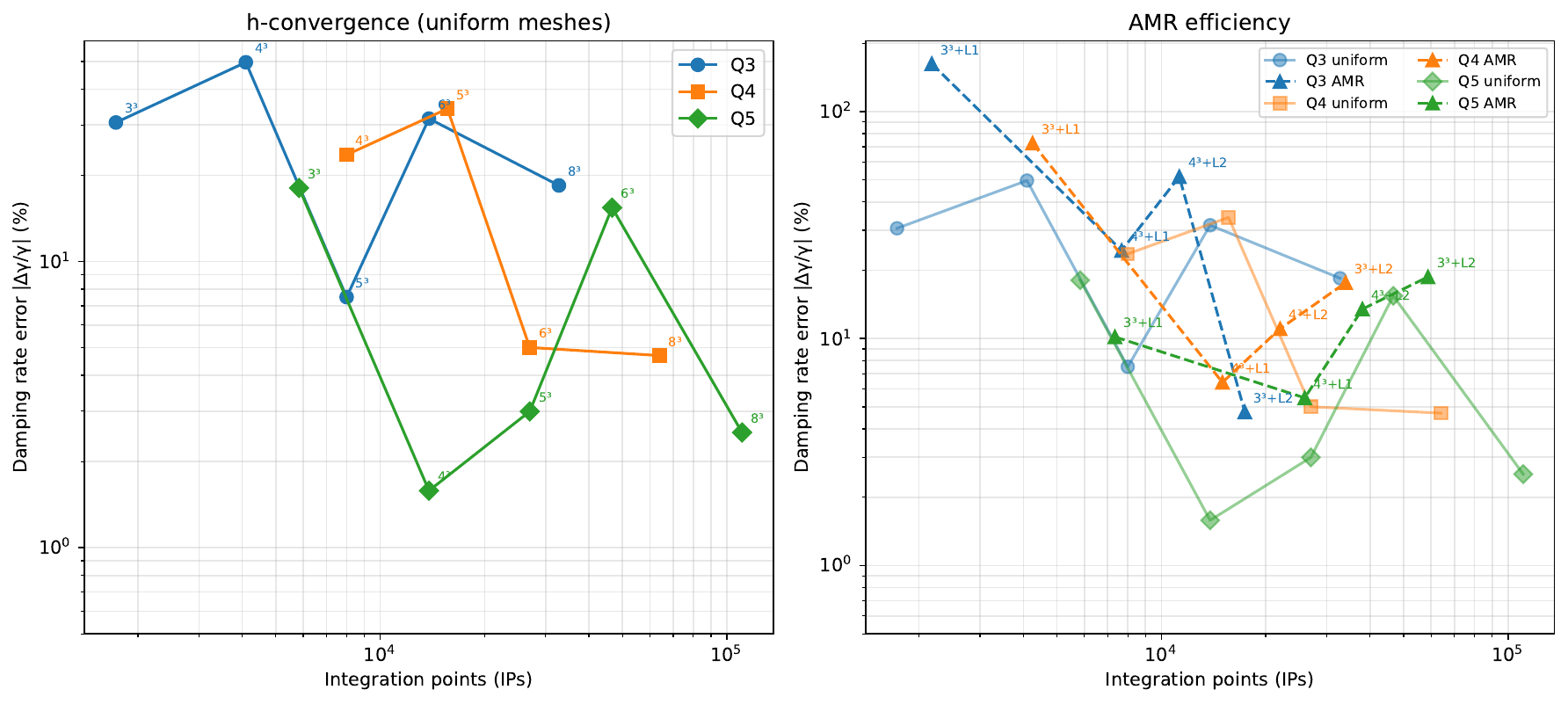}
\caption{Damping rate error vs.\ integration points.
  \emph{Left}: h-convergence on uniform meshes ($L=0$) for $p=3,4,5$;
  each point is labeled with the base grid size $N_b^3$.
  \emph{Right}: AMR efficiency --- uniform meshes (solid, faded) vs.\
  AMR meshes (dashed, triangles) at each polynomial degree.  All runs:
  $R=6$, $N_x=64$, $\Delta t=0.1$, $t_{\rm end}=20$.}
\label{fig:convergence}
\end{figure}

Table~\ref{tab:convergence} reports the fitted damping rate error for
all configurations tested: uniform meshes with $N_b = 3, 4, 5, 6, 8$
and AMR meshes with $N_b = 3, 4$ and $L = 1, 2$ levels of refinement,
for polynomial degrees $p = 3, 4, 5$.
Figure~\ref{fig:convergence} plots the damping rate error versus IP
count on a log-log scale.

\begin{table}[htbp]
\centering
\caption{Damping rate error for uniform and AMR velocity meshes.
  Analytical rate: $\gamma_{\rm th} = -0.1533$ ($k=0.5$).  All runs:
  $R=6$, $N_x=64$, $\Delta t=0.1$, $t_{\rm end}=20$.  Rate error
  $= |(\hat\gamma - \gamma_{\rm th})/\gamma_{\rm th}| \times 100\%$.
  Entries marked ``---'' have fewer than 2 monotonically-decreasing
  peaks, which is insufficient for a reliable envelope fit.
  Upper block: uniform meshes; lower block: AMR meshes.}
\label{tab:convergence}
\smallskip
\begin{tabular}{cccrrrr}
\toprule
$p$ & $N_b^3$ & $L$ & Cells & IPs & \#peaks & Error (\%) \\
\midrule
\multicolumn{7}{c}{\emph{Uniform meshes}} \\
\midrule
  3 & $3^3$ & 0 &   27 &   1{,}728 & 2 & 30.6 \\
  3 & $4^3$ & 0 &   64 &   4{,}096 & 2 & --- \\
  3 & $5^3$ & 0 &  125 &   8{,}000 & 2 &  7.5 \\
  3 & $6^3$ & 0 &  216 &  13{,}824 & 3 & 31.6 \\
  3 & $8^3$ & 0 &  512 &  32{,}768 & 4 & 18.4 \\
\cmidrule{1-7}
  4 & $3^3$ & 0 &   27 &   3{,}375 & 1 & --- \\
  4 & $4^3$ & 0 &   64 &   8{,}000 & 3 & 23.6 \\
  4 & $5^3$ & 0 &  125 &  15{,}625 & 3 & 34.2 \\
  4 & $6^3$ & 0 &  216 &  27{,}000 & 5 &  5.0 \\
  4 & $8^3$ & 0 &  512 &  64{,}000 & 5 &  4.7 \\
\cmidrule{1-7}
  5 & $3^3$ & 0 &   27 &   5{,}832 & 2 & 18.1 \\
  5 & $4^3$ & 0 &   64 &  13{,}824 & 3 &  1.6 \\
  5 & $5^3$ & 0 &  125 &  27{,}000 & 4 &  3.0 \\
  5 & $6^3$ & 0 &  216 &  46{,}656 & 6 & 15.4 \\
  5 & $8^3$ & 0 &  512 & 110{,}592 & 8 &  2.5 \\
\midrule
\multicolumn{7}{c}{\emph{AMR meshes}} \\
\midrule
  3 & $3^3$ & 1 &   34 &   2{,}176 & 2 & --- \\
  3 & $3^3$ & 2 &  272 &  17{,}408 & 4 &  4.7 \\
  3 & $4^3$ & 1 &  120 &   7{,}680 & 5 & 24.4 \\
  3 & $4^3$ & 2 &  176 &  11{,}264 & 6 & --- \\
\cmidrule{1-7}
  4 & $3^3$ & 1 &   34 &   4{,}250 & 3 & --- \\
  4 & $3^3$ & 2 &  272 &  34{,}000 & 6 & 17.6 \\
  4 & $4^3$ & 1 &  120 &  15{,}000 & 5 &  6.4 \\
  4 & $4^3$ & 2 &  176 &  22{,}000 & 5 & 11.0 \\
\cmidrule{1-7}
  5 & $3^3$ & 1 &   34 &   7{,}344 & 4 & 10.2 \\
  5 & $3^3$ & 2 &  272 &  58{,}752 & 7 & 18.7 \\
  5 & $4^3$ & 1 &  120 &  25{,}920 & 7 &  5.5 \\
  5 & $4^3$ & 2 &  176 &  38{,}016 & 8 & 13.5 \\
\bottomrule
\end{tabular}
\end{table}

\paragraph{Pre-asymptotic regime.}
The h-convergence curves in Figure~\ref{fig:convergence} (left) show a
general downward trend with increasing IP count, but the convergence is
non-monotonic: the error oscillates rather than decreasing smoothly.
This indicates that the velocity meshes tested here are in a
\emph{pre-asymptotic regime} where the error is not yet dominated by
the leading-order truncation term.  The Landau damping physics involves
resonant particles near the phase velocity $v_{\rm ph} = \omega_r / k
\approx 2.8$, which falls in the transition region between the
well-resolved Maxwellian core and the coarsely-resolved tail.  At these
moderate resolutions, the interplay between core and tail resolution
produces the observed non-monotonic behavior.  
The AMR strategy in PETSc is optimized with more scientifically relevant tests.
Nevertheless, the best
configurations at each polynomial degree achieve damping rate errors of
$7.5\%$ (Q3, $5^3$), $4.7\%$ (Q4, $8^3$), and $1.6\%$ (Q5, $4^3$),
demonstrating that the SLDG method correctly captures the Landau damping
physics.

\paragraph{AMR efficiency.}
Figure~\ref{fig:convergence} (right) compares AMR and uniform meshes.
The $4^3$+AMR1 configurations are competitive with uniform meshes at
similar IP counts: for example, Q4 $4^3$+AMR1 ($15{,}000$ IPs, $6.4\%$
error) is competitive with Q4 $6^3$ uniform ($27{,}000$ IPs, $5.0\%$ error) with
$1.8\times$ fewer IPs.  However, the $3^3$-based AMR meshes perform
poorly because the $3^3$ base grid has too few coarse cells to resolve
the distribution function in the outer velocity space.  This confirms
that AMR is most effective when the base grid already provides adequate
resolution of the bulk distribution, and the refinement adds resolution
where the gradients are steepest.

\begin{figure}[htbp]
\centering
\includegraphics[width=\textwidth]{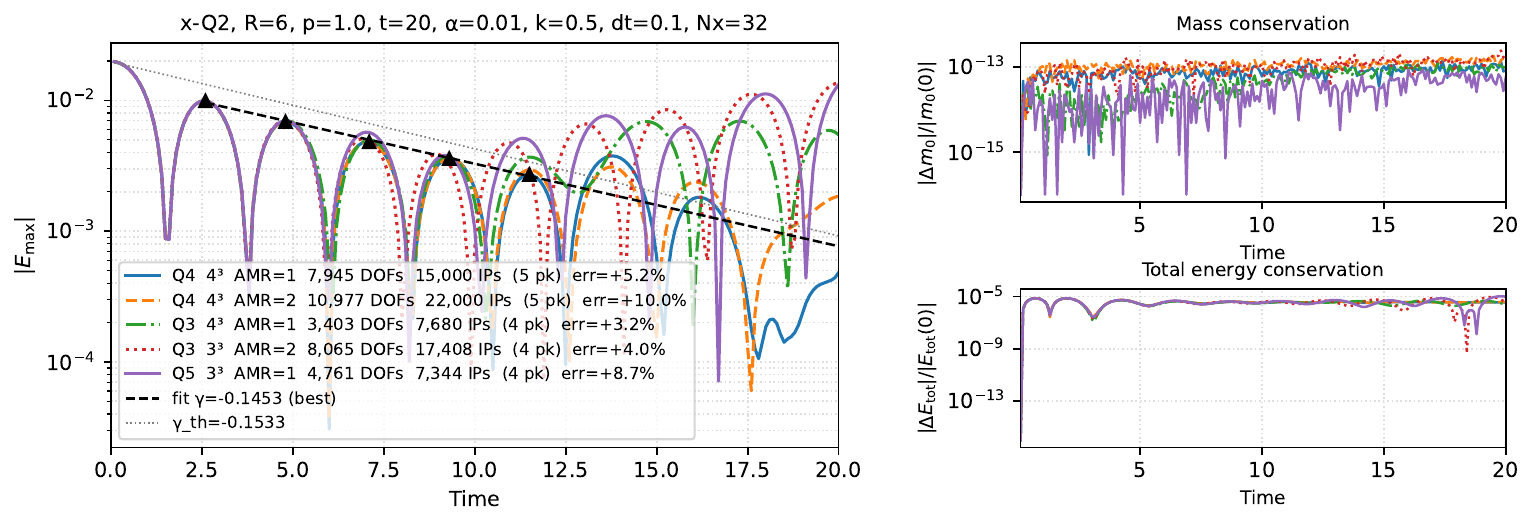}
\caption{\emph{Left}: $E_{\max}$ vs.\ time (semilogy) for selected
  configurations; the black dashed line shows the envelope fit and
  black triangles mark detected peaks for the best configuration
  (Q5, $4^3$ uniform).
  \emph{Right top}: relative mass conservation error
  $|\Delta m_0|/|m_0(0)|$.  \emph{Right bottom}: relative total-energy
  conservation error $|\Delta E_{\rm tot}|/|E_{\rm tot}(0)|$.  All
  runs: $R=6$, $N_x=64$, $\Delta t=0.1$, $t_{\rm end}=20$.}
\label{fig:sweep-best}
\end{figure}

Figure~\ref{fig:sweep-best} shows $E_{\max}$ time series and
conservation diagnostics for representative configurations.

\paragraph{Comparison with PIC.}
For context, Finn et al.\ \cite{FinnLandauDamping2023} solve the same
1X+3V Landau damping benchmark ($k=0.5$, $\alpha=0.01$) using PETSc-PIC
with $10^6$ particles and report the theoretical damping rate
$\gamma = -0.1533$ to within ${\sim}1\%$.  Our best SLDG configuration
(Q5, $4^3$ uniform, $13{,}824$ IPs) achieves $1.6\%$ error.  The SLDG
method uses far fewer degrees of freedom than PIC (${\sim}10^4$ vs.\
$10^6$) and provides exact mass conservation (to machine precision),
whereas PIC conserves mass by construction but suffers from particle
noise that limits the accuracy of higher-order moments.  A detailed
comparison at matched computational cost is left for future work.

\paragraph{Long-time conservation.}
To stress-test conservation over longer integration times, we run the Q3
$4^3$+AMR1 configuration ($120$ cells, $7{,}680$ IPs) for $1{,}000$ time
steps to $t = 100$.  Mass is preserved to a relative error of
$3.4 \times 10^{-13}$ (near machine precision), confirming that the
partition-of-unity property (Section~\ref{sec:amr-conservation}) holds
exactly in exact arithmetic and to machine precision in floating-point.  The total energy
$\mathcal{E} = \frac{1}{2} m_2 + \mathcal{E}_E$ drifts by a relative
error of $1.5 \times 10^{-5}$ over the full run.  Unlike mass, energy
is not an exact invariant of the SLDG discretization; the observed drift
is consistent with the $O(\Delta t^2)$ Strang splitting error
accumulating over $1{,}000$ steps, and remains small and bounded with no
secular growth.

\subsection{Performance}
\label{sec:performance}

Inspecting Figure~\ref{fig:sweep-best}, the configuration with the most
regular $E_{\max}$ decay (most peaks before recurrence) is $p_v = 3$
with a $3^3$ base mesh and 2 AMR levels.
For scaling tests we use the $4^3$ and one level of adaptive refinement (AMR1) configuration ($p_v = 3$, 120
cells, $7{,}680$ IPs), which provides a good balance between accuracy
and problem size.  To obtain a problem size large enough for meaningful
scaling tests, we apply one level of uniform post-refinement: the
command-line option \texttt{-dm\_landau\_amr\_post\_refine 1} uniformly
refines every cell of the AMR mesh once after the adaptive refinement is
complete, increasing the cell count from 120 to 960 without altering the
AMR topology or the fast/slow classification.

\subsection{GPU execution model}
\label{sec:gpu-model}

The solver keeps the distribution function on the GPU throughout the
time-stepping loop.  A persistent Kokkos device View
(\texttt{d\_f\_scratch}) holds the authoritative copy of $f$; the host
PETSc Vec is only updated on demand (e.g., for file output).  All
compute kernels---v-advection pencil sweeps, x-advection SLDG updates,
charge density reduction, and moment diagnostics---execute entirely on
the device via \texttt{Kokkos::parallel\_for} or
\texttt{parallel\_reduce}.  The only host$\leftrightarrow$device
transfers per time step are:
\begin{enumerate}
\item the x-advection ghost exchange: a few ghost cells per MPI rank
  boundary ($\mathrm{sw} \times N_v^{\mathrm{DOF}} \times (p_x+1)$
  doubles per side), and
\item the charge density vector $\rho$ ($N_x^{\mathrm{local}}$ doubles)
  needed by the host-side Poisson solver.
\end{enumerate}
For the Q3 $4^3$+AMR1 configuration with $N_x = 1024$ on 4 GPUs, the
ghost exchange transfers ${\sim}1.7$\,MB per rank boundary per direction
(ghost width $\mathrm{sw}$ cells $\times$ $N_v^{\mathrm{DOF}} \times
(p_x+1)$ doubles per side), and
$\rho$ is ${\sim}2$\,KB.

\paragraph{Persistent device Views.}
All device-side arrays (the distribution function, output buffer,
overlap matrices, pencil metadata, and staging buffers) are allocated
once during setup and reused every step, avoiding per-step
\texttt{cudaMalloc}/\text{cudaFree} overhead.  At the end of the
v-advection step, the input and output Views are swapped via pointer
exchange (\texttt{std::swap}), an $O(1)$ operation that avoids copying
the full distribution function.

\paragraph{Ghost exchange with pinned staging buffers.}
\label{sec:ghost-exchange}
The x-advection ghost exchange uses persistent \emph{pinned host}
staging buffers allocated with
\texttt{Kokkos::View<\ldots, Kokkos::HostSpace>}, which maps to
\texttt{cudaMallocHost} on CUDA builds.  Pinned (page-locked) memory is
directly accessible by the network fabric for RDMA transfers,
eliminating an extra internal copy by the MPI library.  The exchange
proceeds as:
(1)~\texttt{Kokkos::deep\_copy} from device to pinned send buffers,
(2)~\texttt{MPI\_Sendrecv} on pinned buffers, and
(3)~\texttt{Kokkos::deep\_copy} from pinned receive buffers to device.
Interior cells are copied device-to-device.

\subsection{Strong and weak scaling}
\label{sec:strong-scaling}

\begin{figure}[htbp]
\centering
\begin{minipage}[t]{0.48\textwidth}
  \centering
  \includegraphics[width=\textwidth]{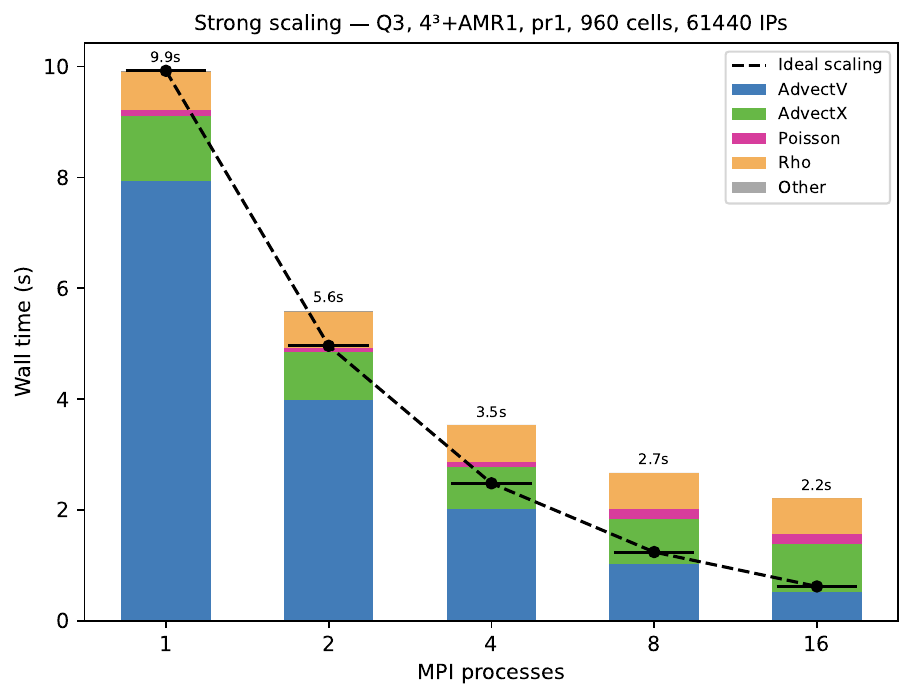}
\end{minipage}
\hfill
\begin{minipage}[t]{0.48\textwidth}
  \centering
  \includegraphics[width=\textwidth]{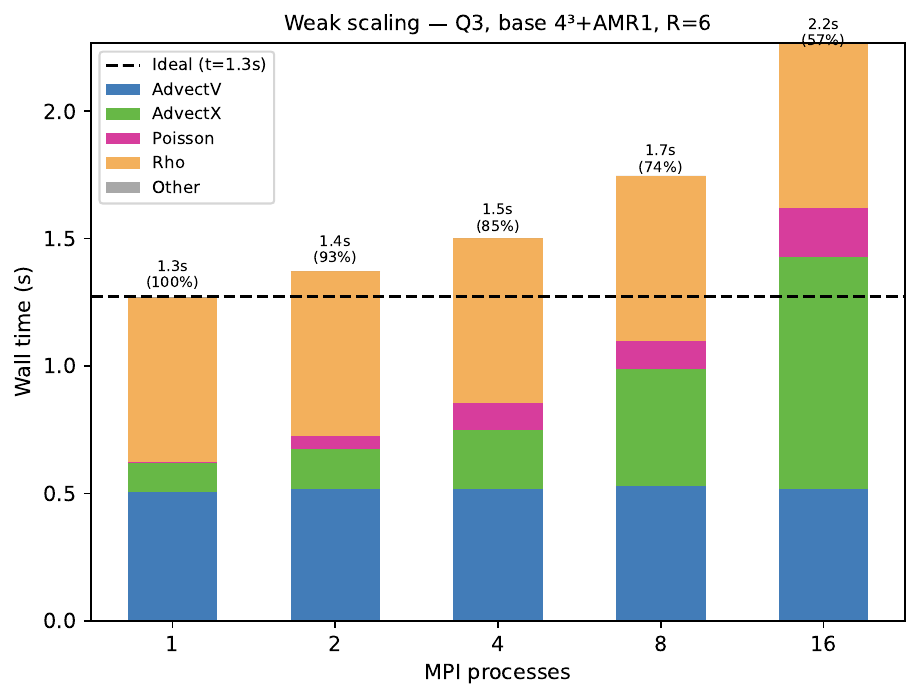}
\end{minipage}
\caption{Scaling of the SLDG solver on Perlmutter GPU nodes (Q3,
  $4^3$+AMR1+post-refine~1, 960 velocity cells, 61{,}440 IPs;
  device-resident $f$ with pinned ghost exchange, see
  Section~\ref{sec:gpu-model}).
  Stacked bar charts show wall time by component (AdvectV, AdvectX,
  Poisson, Rho, Other); the dashed line marks the ideal scaling target.
  \emph{Left} (strong scaling, $N_x=1024$, 20 steps, 1--16 MPI
  processes): efficiency at 16 processes is $28\%$; the charge
  density computation (Rho) is a constant-time bottleneck
  (${\sim}0.65$\,s) that does not scale with process count, and
  AdvectX ghost-exchange communication limits further improvement.
  \emph{Right} (weak scaling, 960 cells per process, $N_x$ scaled
  proportionally, 1--16 MPI processes): efficiency at 16 processes is
  $57\%$, dominated by AdvectX communication overhead and the
  constant Rho cost.}
\label{fig:scaling}
\end{figure}

We measure strong scaling by fixing the problem size at Q3, $4^3$+AMR1
with post-refine level~1 (960 velocity cells, 61{,}440 integration
points) and varying the number of MPI processes from 1 to 16.  Each run
uses 20 time steps with $N_x = 1024$ x-cells.
For weak scaling we keep the per-process velocity work constant at
960 cells ($N_x = 64$ per process) and scale from 1 to 16 MPI processes
by increasing $N_x$ proportionally.
Figure~\ref{fig:scaling} shows the wall-time breakdown by component for
both experiments.

The component breakdown in Figure~\ref{fig:scaling} reveals three
scaling bottlenecks.  At $n_p = 1$, the v-advection kernel (AdvectV)
dominates wall time ($7.9$\,s out of $9.9$\,s total, $80\%$).  AdvectV
scales nearly perfectly: $7.9$\,s $\to$ $0.52$\,s at $n_p = 16$
($15.3\times$ speedup on 16 processes), because each velocity pencil
sweep is embarrassingly parallel across x-cells.  In contrast, AdvectX
scales to only $1.3\times$ at $n_p = 16$ ($1.17$\,s $\to$ $0.87$\,s),
because the ghost exchange grows with the number of MPI rank boundaries
while the per-rank compute (AdvectX\_Apply) shrinks.
A third bottleneck is the
charge density computation (Rho), which remains nearly constant at
${\sim}0.65$\,s regardless of process count.  Detailed profiling shows
that $99.8\%$ of the Rho time is spent in the GPU kernel itself (fence
time), with kernel launch overhead (${\sim}30\,\mu$s) and the
device-to-host copy of the scalar result (${\sim}20\,\mu$s) both
negligible.  The bottleneck is GPU under-occupancy: the current kernel
launches only $N_x^{\mathrm{local}}$ threads (one per x-cell), each
performing a serial reduction over all $N_v^{\mathrm{DOF}} \times
(p_x+1)$ velocity DOFs---far fewer threads than needed to saturate the
A100's compute units.
By $n_p = 8$, AdvectV ($1.0$\,s), AdvectX
($0.83$\,s), and Rho ($0.66$\,s) are all comparable, and further
scaling is limited by the non-scalable Rho and communication overheads.

\section{Conclusions}
\label{sec:conclusions}

We have extended the semi-Lagrangian discontinuous Galerkin (SLDG) method
of Einkemmer \cite{Einkemmer2016} to adaptively refined velocity grids and
to three-dimensional velocity space.  The key algorithmic contributions
are:
\begin{itemize}
\item A \emph{hybrid fast/slow sweep} that uses precomputed per-level
  overlap matrices for conforming cells (identical cost to the uniform-grid
  method) and generalized overlap integrals for nonconforming cells at
  refinement boundaries.

\item A \emph{CSR pencil data structure} that organizes the dimensional
  splitting on AMR meshes, with weighted accumulation for coarse cells
  appearing in multiple pencils.

\item Extension to \emph{1X+3V} using tensor-product DG elements on
  hexahedral cells with an efficient pack/scatter strategy based on a
  precomputed tensor-product permutation map.
\end{itemize}
The solver preserves the key properties of the original SLDG method:
conservation (mass preserved to machine precision with periodic boundary
conditions), a compact two-cell stencil, and low numerical diffusion.  On
the Landau damping benchmark, the Q4 $4^3$+AMR1 mesh ($15{,}000$ IPs)
achieves comparable damping rate accuracy to the Q4 $6^3$ uniform mesh
($27{,}000$ IPs) with $1.8\times$ fewer integration points.  Total energy
drifts by $O(10^{-5})$ over $1{,}000$ time steps ($t=100$); mass is
conserved to machine precision.  On Perlmutter A100 GPUs, the Kokkos-based
implementation with pinned host staging buffers achieves $28\%$ strong
scaling efficiency and $57\%$ weak scaling efficiency at 16 processes
for the Q3 $4^3$+AMR1 configuration with 61{,}440 integration points;
the primary scaling bottlenecks are the charge density computation (Rho),
which suffers from GPU under-occupancy and does not benefit from
additional processes, and the x-advection ghost exchange, whose MPI
communication cost grows with the number of rank boundaries.

The solver is implemented in PETSc and reuses the velocity mesh
infrastructure of PETSc's Landau collision operator
\cite{AdamsLandau2025}, enabling future coupling for full
Vlasov--Maxwell--Landau simulations.  Directions for future work include:
\begin{itemize}
\item Extension to multiple spatial dimensions (3X+3V).
\item Coupling with the Landau collision operator for collisional plasma
  simulations.
\item \emph{Charge density kernel optimization.}  The current Rho kernel
  launches one GPU thread per x-cell, leaving the A100 severely
  under-occupied.  Replacing the flat \texttt{parallel\_for} with a
  \texttt{TeamPolicy} reduction (one team per x-cell, with team threads
  cooperatively reducing over velocity DOFs) would increase GPU occupancy
  by ${\sim}250\times$ and is expected to reduce the Rho cost by an order
  of magnitude.
\item \emph{Fusing consecutive x-advection half-steps.}  In the Strang
  splitting, the trailing $\Delta t/2$ x-advection of step~$n$ and the
  leading $\Delta t/2$ of step~$n{+}1$ use the same precomputed overlap
  matrices (since the x-advection speed $v_x$ is independent of the
  electric field).  Fusing these into a single $\Delta t$ step would
  eliminate one ghost exchange per time step, halving the communication
  cost.
\item \emph{Non-blocking MPI ghost exchange.}  Replacing the current
  blocking \texttt{MPI\_Sendrecv} calls with non-blocking
  \texttt{MPI\_Isend}/\texttt{MPI\_Irecv} would overlap the two exchange
  directions and could be combined with an interior/boundary kernel split
  to overlap communication with computation.
\end{itemize}

\paragraph{Reproducibility.}
Data, plotting scripts, and run commands for all numerical results in this
paper are available at \url{https://gitlab.com/markadams4/sldg-amr}.

\paragraph{Use of AI tools.}
Claude (Anthropic) was used in the development of this work, including
code development and manuscript preparation.

\section*{Acknowledgments}

This material is based upon work supported by the U.S.\ Department of
Energy, Office of Science, Office of Advanced Scientific Computing
Research, Scientific Discovery through Advanced Computing (SciDAC) Program
through the FASTMath Institute, under contract number DE-AC02-05CH11231 at
Lawrence Berkeley National Laboratory.

\bibliographystyle{siam}
\bibliography{references}

\end{document}